\documentclass[11pt]{amsart}

\usepackage{url}

\newcommand{\R}{\mathbb R}
\newcommand{\Z}{\mathbb Z}

\newcommand{\cand}[1]{\mathsf{#1}}

\newcommand{\MP}{{\mathcal P}}
\newcommand{\MS}{{\mathcal S}}

\DeclareMathOperator{\card}{card}

\DeclareMathOperator{\Prob}{Prob}

\DeclareMathOperator{\relvol}{relvol}
\DeclareMathOperator{\lt}{lt}

\title[Exploiting Polyhedral Symmetries in Social Choice]{Exploiting Polyhedral Symmetries\\ in Social Choice} 
\author{Achill Sch\"urmann} 
\address{%
Institute of Mathematics\\
University of Rostock\\
18051 Rostock\\
Germany}
\email{achill.schuermann@uni-rostock.de}

\subjclass[2010]{91B, 52B}

\keywords{Social Choice,  Impartial Anonymous Culture, Condorcet, plurality voting, symmetric polyhedra, weighted Ehrhart theory}

\begin{document}

\begin{abstract}
A large amount of literature in social choice theory deals
with quantifying the probability of certain election outcomes.
One way of computing the probability of a specific voting situation
under the Impartial Anonymous Culture assumption
is via counting integral points in polyhedra.
Here, Ehrhart theory can help, but unfortunately the dimension
and complexity of the involved polyhedra grows rapidly with the
number of candidates.
However, if we exploit available polyhedral symmetries,
some computations become possible that previously
were infeasible.
We show this in three well known examples:
Condorcet's paradox, Condorcet efficiency of plurality voting and in Plurality voting vs Plurality Runoff.
\end{abstract}

\maketitle

\tableofcontents

\section{Introduction}

In social choice theory, a vast amount of literature deals 
with quantifying the probability of certain election outcomes.
This is in particular the case for so-called ``voting paradoxes''
that are known to be unavoidable since the famous 
Impossibility Theorem of Arrow~\cite{arrow-1951}
(see \cite{tp-2008} for a popular exposition).
Under the Impartial Anonymous Culture (IAC) assumption,
the probability for such an event can be computed by
counting integral solutions to a system of linear inequalities,
associated to the specific voting event of interest (see for
example~\cite{gl-2011}).
There exists a rich mathematical theory 
going back to works of Ehrhart~\cite{ehrhart-1967} in the 1960s 
that helps to deal with such problems.
We refer to \cite{br-2007} and \cite{barvinok-2008} for an introduction.
The connection to Social Choice Theory was recently discovered 
by Lepelley et al. \cite{lls-2008}  and Wilson and Pritchard \cite{wp-2007}. 
A few years earlier a similar theory had been described specifically
for the social choice context by Huang and Chua~\cite{hc-2000} (see
also \cite{gehrlein-2002}).
Based on Barvinok's algorithm \cite{barvinok-1994} there now exists 
specialized mathematical software 
for performing previously cumbersome or practically impossible computations.
The first available program was {\tt LattE}, 
with its newest version  {\tt LattE integrale} (see~\cite{latte});
an alternative implementation of Barvinok's algorithm is available 
through {\tt barvinok} (see~\cite{barvinok}), which is also usable within the 
{\tt polymake} framework (see~\cite{polymake}).

The purpose of this note is to shed some light on the 
possibilities for social choice computations that arise 
through the use of Ehrhart theory and 
weighted generalizations of it~(see~\cite{bbdmv-2010}). 
We in particular show how
symmetry of linear systems characterizing certain voting events, 
can be used to obtain new results.
As examples, we consider three well studied 
voting situations with four candidates:
{\em Condorcet's paradox},
the {\em Condorcet efficiency of plurality voting} and different 
outcomes in {\em Plurality vs Plurality Runoff}.

In Section~\ref{sec:three-cand} we review some linear models
for voting events and introduce some of the used notation.
In Section~\ref{sec:likelihood} we sketch how counting integral points
in polyhedra and Ehrhart's theory can be used to compute probabilities
for voting outcomes. In Section~\ref{sec:reduce-dim} we
show how the complexity of computations can be reduced by 
using a symmetry reduced, lower dimensional reformulation. 
We in particular show how to use
integration to obtain exact
values for the limiting probability of voting outcomes
when the number of voters tends to infinity.
As examples, we obtain previously unknown, exact values for
election events with four candidates.

%
%

\section{Linear systems describing voting situations}

\label{sec:three-cand}

\subsection*{Notation}

For the start we look at three candidate elections, as everything that
will follow can best be motivated and explained in smaller examples.
Assume there are $n$~voters, with $n \geq 2$, and each of them has a 
complete linear (strict) preference order on the three candidates $\cand{a},\cand{b},\cand{c}$.
We subdivide the voters into six groups 
\begin{equation} \label{eqn:six-voters}
\left(n_{\mathsf{ab}},n_{\mathsf{ac}},n_{\mathsf{ba}},n_{\mathsf{bc}},n_{\mathsf{ca}},n_{\mathsf{cb}}\right)
,
\end{equation}
according to their six possible preference orders:
$$
\cand{a}\cand{b}\cand{c} (n_{\cand{a}\cand{b}}) \qquad 
\cand{a}\cand{c}\cand{b} (n_{\cand{a}\cand{c}}) \qquad
\cand{b}\cand{a}\cand{c} (n_{\cand{b}\cand{a}}) \qquad 
\cand{b}\cand{c}\cand{a} (n_{\cand{b}\cand{c}}) \qquad
\cand{c}\cand{a}\cand{b} (n_{\cand{c}\cand{a}}) \qquad 
\cand{c}\cand{b}\cand{a} (n_{\cand{c}\cand{b}}) 
$$
For example, there are $n_{\cand{a}\cand{b}}$ voters that prefer
$\cand{a}$ over $\cand{b}$ and $\cand{b}$ over $\cand{c}$.
We omit the last preference in the index, as it is determined 
once we know the others. This type of indexing will
show to be useful when we reduce the number of variables
in Section~\ref{sec:reduce-dim}.

The tuple~\eqref{eqn:six-voters} is referred to as a {\em voting situation}.
In an election with 
\begin{equation}  \label{eqn:n-equation}
n=n_{\mathsf{ab}} + n_{\mathsf{ac}} +
n_{\mathsf{ba}} + n_{\mathsf{bc}} + n_{\mathsf{ca}} + n_{\mathsf{cb}}
\end{equation}
voters, there are $n+5 \choose 5$ possible voting situations.
We make the simplifying {\em Impartial Anonymous Culture (IAC)
  assumption} that each of these voting situations is equally 
likely to occur.

\subsection*{Condorcet's Paradox}

Maybe the most famous voting paradox 
goes back to the Marquis de Condorcet (1743--1793).
He observed that in an election with three or more candidates, it is
possible that pairwise comparison of candidates can lead to an
intransitive collective choice. For instance, candidate~$\cand{a}$ could be
preferred over candidate~$\cand{b}$, $\cand{b}$ could be preferred
over candidate $\cand{c}$ and $\cand{c}$ could be preferred over
candidate $\cand{a}$.
In this case there is no 
{\em Condorcet winner}, that is, someone who beats every other candidate by
pairwise comparison.

The condition that candidate~$\cand{a}$ is a Condorcet winner 
can be described via two linear constraints:

\noindent
\begin{minipage}[h]{0.7\textwidth}
\begin{eqnarray}  \label{eq:a-beats-b}
n_{\mathsf{ab}} + n_{\mathsf{ac}} + n_{\mathsf{ca}}  & > &
n_{\mathsf{ba}} + n_{\mathsf{bc}} + n_{\mathsf{cb}}
\\
\label{eq:a-beats-c}
n_{\mathsf{ab}} + n_{\mathsf{ac}} + n_{\mathsf{ba}} & > &
n_{\mathsf{ca}} + n_{\mathsf{bc}} + n_{\mathsf{cb}} 
\end{eqnarray}
\end{minipage}
\hfill
\begin{minipage}[h]{0.3\textwidth}
\begin{eqnarray*}
 \mbox{( $\cand{a}$ beats $\cand{b}$ )} \\
 \mbox{ ( $\cand{a}$ beats $\cand{c}$ )}
\end{eqnarray*}
\end{minipage}

\medskip

The probability of candidate $\cand{a}$ being a Condorcet winner in an
election with $n$ voters can be expressed as the quotient 
$$
\Prob(n) \; = \;
\frac{\card\left\{\left(n_{\mathsf{ab}},\ldots,n_{\mathsf{cb}}\right)
  \in \Z^6_{\geq 0}  \mbox{ satisfying \eqref{eqn:n-equation},
    \eqref{eq:a-beats-b}, \eqref{eq:a-beats-c} } \right\}}{{n+5
  \choose 5}}
.
$$
The denominator is a polynomial of degree $5$ in $n$.
It had been observed by Fishburn and Gehrlein \cite{gf-1976} (cf.~\cite{bb-1983}) that the
numerator shows a similar behavior: Restricting to even or odd $n$ it
can be expressed as a degree $5$ polynomial in $n$. The leading
coefficient of both polynomials is the same and we approach the same
probability for large elections (as $n$ tends to infinity). This
{\em limiting probability} is known to be
$$
\lim_{n\to\infty}\Prob(n) = \frac{5}{16}
.
$$

Having the probability for candidate $\cand{a}$ being a Condorcet
winner, we obtain the probability for a Condorcet paradox (no
Condorcet winner exists) as
$1- 3\cdot \Prob(n)$ with an exact limiting probability of $\tfrac{1}{16}$.

\medskip

In a similar way we can determine probabilities for other voting events.

\subsection*{Condorcet efficiency of Plurality voting}

If there is a Condorcet winner, there is good reason 
to consider him to be the voter's choice. However, many common 
voting rules do not always choose the Condorcet winner even
if one exists. This is in particular the case for the widely 
used plurality voting, where the candidate 
with a majority of first preferences is elected.

The condition that candidate~$\cand{a}$ is a Condorcet winner
but  candidate~$\cand{b}$ is the plurality winner can be expressed
by the two inequalities~\eqref{eq:a-beats-b} and~\eqref{eq:a-beats-c}, together with the two additional inequalities

\noindent
\begin{minipage}[h]{0.5\textwidth}
\begin{eqnarray} \label{eq:b-wins-plurality-over-a}
n_{\mathsf{ba}} + n_{\mathsf{bc}}  & > &
n_{\mathsf{ab}} + n_{\mathsf{ac}} 
\\ 
\label{eq:b-wins-plurality-over-c}
n_{\mathsf{ba}} + n_{\mathsf{bc}} & > &
n_{\mathsf{ca}} +  n_{\mathsf{cb}}
\end{eqnarray}
\end{minipage}
\hfill
\begin{minipage}[h]{0.5\textwidth}
\begin{eqnarray*}
 \mbox{( $\cand{b}$ wins plurality over $\cand{a}$ )} \\
 \mbox{ ( $\cand{b}$ wins plurality over $\cand{c}$ )}
\end{eqnarray*}
\end{minipage}

\medskip

The {\em Condorcet efficiency} of a voting rule is the conditional probability
that a Condorcet winner is elected if one exists. As there are $3\cdot
2$ possibilities for choosing a Condorcet winner
and another plurality winner, we obtain 
$$
\Prob(n) \; = \;
\frac{6 \cdot \card\left\{\left(n_{\mathsf{ab}},\ldots,n_{\mathsf{cb}}\right)
  \in \Z^6_{\geq 0}  \mbox{ satisfying \eqref{eqn:n-equation},
    \eqref{eq:a-beats-b}, \eqref{eq:a-beats-c},
    \eqref{eq:b-wins-plurality-over-a}, \eqref{eq:b-wins-plurality-over-c}} \right\}}
   {3 \cdot \card\left\{\left(n_{\mathsf{ab}},\ldots,n_{\mathsf{cb}}\right)
   \in \Z^6_{\geq 0}  \mbox{ satisfying \eqref{eqn:n-equation},
     \eqref{eq:a-beats-b}, \eqref{eq:a-beats-c} } \right\}}
$$
for the likelihood of a Condorcet winner being a plurality loser.
Again, depending on $n$ being odd or even, one obtains polynomials
in~$n$ in the denominator and the numerator (see~\cite{gehrlein-1982}).
The exact value of the limit $\lim_{n\to\infty} \Prob(n)$ is $16/135$. 
Therefore, the Condorcet efficiency of plurality voting with three candidates
is $119/135 = 88.\overline{148}\%$.

\subsection*{Plurality vs Plurality Runoff}

Plurality Runoff voting is a common practice to overcome some of these ``problems'' of Plurality voting.
It is used in many presidential elections, for example in
France. After a first round of plurality voting in which none of the
candidates has achieved more than $50\%$ of the votes,
the first two candidates compete in a second runoff round.

The condition that candidate~$\cand{b}$ is the plurality winner,
but candidate~$\cand{a}$ wins the second round of Plurality Runoff
can be expressed by the two
inequalities~\eqref{eq:b-wins-plurality-over-a} and

\noindent
\begin{minipage}[h]{0.5\textwidth}
\begin{equation} \label{eq:a-wins-plurality-over-c}
n_{\mathsf{ab}} + n_{\mathsf{ac}}   \; > \; 
n_{\mathsf{ca}} + n_{\mathsf{cb}}
, 
\end{equation}
\end{minipage}
\hfill
\begin{minipage}[h]{0.5\textwidth}
$$
 \mbox{( $\cand{a}$ wins plurality over $\cand{c}$ )} \\
$$
\end{minipage}

\medskip

\noindent
together with the linear condition \eqref{eq:a-beats-b} that~$\cand{a}$ beats $\cand{b}$ in a pairwise comparison.
The probability that another candidate is chosen in the second round 
as the number of voters tends to infinity is known to be 
$71/576 = 12.3263\overline{8}\%$ (see \cite{lls-2008}).

\subsection*{Four and more candidates}

\label{sec:four-cand}

Having $m$ candidates we can set up similar linear systems in
$m!$ variables. 
For example, in an election with four candidates
$\cand{a}, \cand{b}, \cand{c}, \cand{d}$ we use the $24$-dimensional
vector $x^t=( n_{\mathsf{abc}},\ldots,  n_{\mathsf{dcb}} )$. 
Here, indices are taken in lexicographical order. 
The condition
that $\cand{a}$ is a Condorcet winner is described by the three
inequalities that imply $\cand{a}$ beats $\cand{b}$, $\cand{a}$ beats
$\cand{c}$
and $\cand{a}$ beats $\cand{d}$ in a pairwise comparison.
As linear systems with $24$~variables become hard to grasp,
it is convenient to use matrices for their description.
We are interested in all non-negative integral (column) vectors $x$ 
satisfying the matrix inequality $Ax>0$ for the matrix 
$A\in \Z^{3\times 24}$ with entries

\medskip

\noindent
\begin{minipage}{0.5cm}
\begin{equation} \label{eqn:big-matrix}
\phantom{bla}
\end{equation}
\end{minipage}
\hfill
\begin{minipage}{12cm}
\begin{footnotesize}
\begin{verbatim}
 1  1  1  1  1  1 -1 -1 -1 -1 -1 -1  1  1 -1 -1  1 -1  1  1 -1 -1  1 -1
 1  1  1  1  1  1  1  1 -1 -1  1 -1 -1 -1 -1 -1 -1 -1  1  1  1 -1 -1 -1
 1  1  1  1  1  1  1  1  1 -1 -1 -1  1  1  1 -1 -1 -1 -1 -1 -1 -1 -1 -1
\end{verbatim}
\end{footnotesize}
 \end{minipage}

\section{Likelihood of voting situations and Ehrhart's theory}

\label{sec:likelihood}

\subsection*{Integral points in polyhedral cones}
In order to deal with an arbitrary number of candidates, 
let us put the example above in a slightly more general context.
In any of the three voting examples, the voting situations
of interest lie in a {\em polyhedral cone}, that is, in a
set $\MP$ of points in $\R^d$ (with $d=6$ or $d=24$ in case of three
or four candidate elections)
satisfying a finite number of homogeneous linear inequalities.
In addition to the strict inequalities which are different in each
of the examples, the condition that the variables $n_i$ are
non-negative can be expressed by the homogeneous linear inequalities $n_i\geq 0$.

Let $\MP, \MS \subset \R^d$ denote two $d$-dimensional {\em polyhedral
  cones}, each defined by some homogeneous linear (possibly strict) inequalities.
We may assume that $\MP$ is contained in $\MS$ and that both
polyhedral cones are contained in the orthant $\R_{\geq 0}^d$. 
If we are interested in elections with $n$~voters, we consider the 
voting situations (integral vectors) in the intersection of~$\MP$ and~$\MS$ with
the affine subspace 
$$
L^d_n =
\left\{ 
(n_1,\ldots, n_d) \in \R^d 
\; | \;
\sum_{i=1}^d n_i = n
\right\}
.
$$
The {\em expected frequency} 
of voting situations being in $\MP$ among voting situations in $\MS$
is then expressed by
\begin{equation} \label{eqn:abstract-prob}
\Prob(n)
=
\frac{\card \left(\MP \cap L^d_n \cap \Z^d \right)}{\card \left(\MS \cap L^d_n \cap \Z^d\right)}
.
\end{equation}

When estimating the probability of candidate $\cand{a}$ being a
Condorcet winner for instance, the homogeneous polyhedral cone $\MS$ is simply the
non-negative orthant $\R_{\geq 0}^d$ described by the linear
inequalities $n_i\geq 0$. In that case the denominator is known to be equal to 
$$
\binom{n+d-1}{d-1} 
.
$$
This is a polynomial in $n$ of degree $d-1$ 
(the dimension of $L^d_n \cap \MS$).

\subsection*{Ehrhart theory}

By Ehrhart's theory, 
the number of integral solutions in a polyhedral cone intersected with~$L^d_n$
can be expressed by a
{\em quasi-polynomial} in~$n$. 
Roughly speaking, a quasi-polynomial is simply a
finite collection $p_1(n), \ldots , p_k(n)$ of polynomials, such that
the number of voting situations is given by 
$p_i(n)$ if $ \; i \equiv n \mod k$.

The degree of the polynomial is equal to the dimension of the
polyhedral cone intersected with~$L^d_n$.
In the voting events considered here this dimension is always equal to $d-1$.
So in the examples with three candidates their degree is always~$5$.
The number~$k$ of different polynomials depends on the linear
inequalities involved.
For the Condorcet paradox 
we have $k=2$ polynomials $p_1(n)$ and $p_2(n)$, where
$p_1(n)$ gives the answer for odd~$n$ ($1 \equiv n \mod 2$)
and $p_2(n)$ gives the answer for even~$n$ ($0 \equiv  2 \equiv n \mod
2$).
For Condorcet efficiency we have $k=6$ (see~\cite{gehrlein-2002})
and for Plurality vs Plurality Runoff we have $k=12$ (see~\cite{lls-2008}).

Given a polyhedral cone $\MP$, the  quasi-polynomial 
$q(n)=\card \left(\MP \cap L^d_n \cap \Z^d \right)$ can be explicitly computed 
using software packages like {\tt LattE integrale} \cite{latte-link} or {\tt
  barvinok} \cite{barvinok-link}. The result for the polyhedral cone~$\MP$
describing candidate~$\cand{a}$ as the Condorcet winner  
could look like

\medskip

\begin{center}
\begin{minipage}{8cm}
\begin{verbatim}
   1/384 * n^5
 + ( 1/64 * { 1/2 * n } + 1/32 ) * n^4 
 + ( 17/96 * { 1/2 * n } + 13/96 ) * n^3
 + ( 23/32 * { 1/2 * n } + 1/4 ) * n^2
 + ( 233/192 * { 1/2 * n } + 1/6 ) * n
 + ( 45/64 * { 1/2 * n } + 0 )
\end{verbatim}
\end{minipage}
\end{center}

\medskip

The curly brackets $\{ \cdots \}$ mean the fractional part of the
enclosed number, allowing to write the quasi-polynomial in a closed form.
In this example we get different polynomials for odd and even $n$. 
Note that the leading coefficient (the coefficient of $n^5$) is in
both cases the same. 
By Ehrhart's theory this is always the case, as it is equal to the
{\em relative volume} of the polyhedron $\MP \cap L^d_1$. 
That is, it is equal to a $\sqrt{d}$-multiple of the standard Lebesgue measure on
the affine space $L^d_1$.
The measure is normalized so that the space contains one integral
point per unit volume.

One technical obstacle using software like {\tt LattE integrale} or {\tt barvinok}
is the use of polyhedral cones described by a mixture of strict and
non-strict inequalities. As the software assumes the input to have only 
non-strict inequalities or equality conditions, one has to avoid the use of strict inequalities. A simple way to achieve this is the replacement of 
strict inequalities $x>0$ by non-strict ones $x\geq 1$, in case $x$ is known to be integral. For instance, if $x$ is a linear expression with integer coefficients,
and if we are interested in integral solutions as in our examples, this is a possible reformulation.

%


Altogether, by obtaining quasi-polynomials for numerator and denominator in \eqref{eqn:abstract-prob}
we get an explicit formula for $\Prob(n)$ via Erhart's theory.

\subsection*{Limiting probabilities via integration}

If we want to compute the exact value of 
$\lim_{n\to \infty} \Prob(n)$ as $n$ tends to infinity, we can use volume computations 
without using Ehrhart's theory. As mentioned above, the leading
coefficients of denominator and numerator correspond to the 
relative volumes of the sets $\MP\cap L_1$ and $\MS\cap L_1$:
$$
\lim_{n\to \infty} \Prob(n)
\; = \;
\lim_{n\to \infty} 
\frac{\card \left( \MP \cap L^d_1 \cap (\Z/n)^d \right)}
        {\card \left( \MS \cap L^d_1 \cap (\Z/n)^d \right)}
\; = \;
\frac{\relvol \left(\MP \cap L^d_1 \right)}{\relvol \left(\MS \cap L^d_1\right)}
$$

In fact, as long as we use the same measure to evaluate the numerator
and the denominator, it does not matter what multiple of the standard
Lebesgue measure we use to compute volume on the affine space $L^d_1$.
The exact relative volume can be computed using {\tt LattE integrale}.
Alternatives are for example {\tt Normaliz} (see~\cite{normaliz-link}) or {\tt vinci} (see~\cite{vinci}).
Exact computations can be quite involved in higher dimensions
(cf. \cite{df-1988}). 
In such cases it is sometimes only possible to compute an approximation, using 
{\em Monte Carlo methods} for instance.

\section{Reducing the dimension by exploiting polyhedral symmetries}

\label{sec:reduce-dim}

In many models the involved linear systems and polyhedra 
are quite symmetric. In particular, permutations of variables may lead to
equivalent linear systems describing the same polyhedron. 
Such symmetries are often visible in smaller examples and can
automatically be determined for larger problems, for instance by our software {\tt SymPol} (see~\cite{sympol}). 
In the three examples described in Section~\ref{sec:three-cand}, we
can exploit such symmetries to reduce the complexity of computations.

\subsection*{Condorcet's paradox}

In case of $\cand{a}$ being a Condorcet winner in a three candidate
election, the variables $n_{\mathsf{ab}}$ and $n_{\mathsf{ac}}$ 
occur pairwise (as $n_{\mathsf{ac}} + n_{\mathsf{ab}}$)  in inequalities
\eqref{eq:a-beats-b}, \eqref{eq:a-beats-c} and in equation \eqref{eqn:n-equation}.
The same is true for $n_{\mathsf{bc}}$ and $n_{\mathsf{cb}}$.
By introducing new variables $n_{\mathsf{a}} = n_{\mathsf{ac}} + n_{\mathsf{ab}}$
and  $n_{\mathsf{\ast a}} = n_{\mathsf{bc}} + n_{\mathsf{cb}}$
we can reduce the dimension of the linear system to only four variables:
\begin{eqnarray*}
n_{\mathsf{a}}  + n_{\mathsf{ca}} - n_{\mathsf{\ast a}} - n_{\mathsf{ba}}  & > & 0  \\  
n_{\mathsf{a}}  + n_{\mathsf{ba}} - n_{\mathsf{\ast a}} - n_{\mathsf{ca}}  & > & 0  \\ 
n_{\mathsf{a}}  + n_{\mathsf{ca}} + n_{\mathsf{\ast a}} + n_{\mathsf{ba}} & = & n \\
n_{\mathsf{a}}, n_{\mathsf{\ast a}}, n_{\mathsf{ba}}, n_{\mathsf{ca}}  & \geq & 0 
.
\end{eqnarray*}
The index $\mathsf{a}$ indicates that we group all variables which carry
candidate $\cand{a}$ as their first preference and index
$\mathsf{\ast a}$ stands for grouping of all variables with
candidate $\cand{a}$ ranked last. In the reduced linear system
each $4$-tuple $(n_{\mathsf{a}},n_{\mathsf{\ast a}} , n_{\mathsf{ba}} , n_{\mathsf{ca}})$
represents several voting situations, previously described by $6$-tuples. For
$n_{\mathsf{a}}$ we have $(n_{\mathsf{a}} + 1)$ different possibilities of
non-negative integral tuples $(n_{\mathsf{ac}},n_{\mathsf{ab}})$.
Similar is true for $n_{\mathsf{\ast a}}$. Together we have 
$$
(n_{\mathsf{a}} + 1) (n_{\mathsf{\ast a}} + 1)
$$
voting situations with three candidates represented by each non-negative integral vector $(n_{\mathsf{a}},n_{\mathsf{\ast a}} , n_{\mathsf{ba}} , n_{\mathsf{ca}})$.

In the four candidate case it is possible to obtain a similar
reformulation by grouping among $24$~variables.
We introduce a new variable for sets of variables having same coefficients in
the linear system. Having a matrix
representation as in~\eqref{eqn:big-matrix}, this kind of special symmetry in the linear system is easy to find by identifying equal columns. Introducing a new variable 
for each set of equal columns we get

\begin{minipage}{0.5cm}
\begin{equation} \label{eqn:a-condorcet-winner-among-four}
\phantom{bla}
\end{equation}
\end{minipage}
\hfill
\begin{minipage}{12cm}
\begin{eqnarray*}
n_{\mathsf{a}}  - n_{\mathsf{ba}}  + n_{\mathsf{ca}}   + n_{\mathsf{da}}  + n_{\mathsf{\ast ab}}  - n_{\mathsf{\ast ac}}  - n_{\mathsf{\ast ad}}  - n_{\mathsf{\ast a}} & > & 0  \\ 
n_{\mathsf{a}}  + n_{\mathsf{ba}}  - n_{\mathsf{ca}}   + n_{\mathsf{da}}  - n_{\mathsf{\ast ab}}  + n_{\mathsf{\ast ac}}  - n_{\mathsf{\ast ad}}  - n_{\mathsf{\ast a}} & > & 0  \\ 
n_{\mathsf{a}}  + n_{\mathsf{ba}}  + n_{\mathsf{ca}}   - n_{\mathsf{da}}  - n_{\mathsf{\ast ab}}  - n_{\mathsf{\ast ac}}  + n_{\mathsf{\ast ad}}  - n_{\mathsf{\ast a}} & > & 0 
\end{eqnarray*}
\end{minipage}

\medskip

These three inequalities describe voting situations in which
candidate~$\cand{a}$ beats candidates $\cand{b}$, $\cand{c}$ and
$\cand{d}$ each in a pairwise comparison.
As in all of our examples, we
additionally have the condition that the involved variables add up to
$n$ and that all of them are non-negative.

As before, 
the used indices of variables reflect which voter preferences are
grouped. As in the three candidate case, $n_{\mathsf{a}}$ and $n_{\mathsf{\ast a}}$ denote
the number of voters with candidate $\cand{a}$ being their first and
last preference respectively. Similarly, ${\mathsf{xy}}$ and
${\mathsf{\ast yx}}$ in the index indicate that voters with preference order
starting with $\cand{x}$, $\cand{y}$ 
and ending with $\cand{y}$, $\cand{x}$ have been combined.

Using our software {\tt SymPol} \cite{sympol-link}
one easily checks that the original system with $24$~variables 
has a symmetry group of order~$199065600$.
The new reduced system with $8$ variables, obtained through the
above grouping of variables, turns out to have a symmetry group 
of order~$6$ only. So most of the symmetry in the original system is of the simple form that is detectable through equal columns in a matrix representation.
The remaining $6$-fold symmetry comes from the 
possibility to arbitrarily permute the variables 
$n_{\mathsf{ba}}, n_{\mathsf{ca}}, n_{\mathsf{da}}$ when at the same
time the variables 
$n_{\mathsf{\ast ab}}, n_{\mathsf{\ast ac}}, n_{\mathsf{\ast ad}}$ are
permuted accordingly.
This symmetry is due to the fact that 
candidates $\cand{b}$, $\cand{c}$ and $\cand{d}$ are equally treated 
in the linear system~\eqref{eqn:a-condorcet-winner-among-four}.
The two new variables $n_{\mathsf{a}}$ and $n_{\mathsf{\ast a}}$ each 
combine six of the former variables. The other six new variables
each combine two former ones.

\subsection*{Weighted counting}

In general, 
if we group more than two variables, say if we 
substitute the sum of $k$ variables $n_1 + \ldots + n_k$ by a new
variable $N$, we have to include a factor of 
$$
\binom{N + k-1}{k-1}
$$
when counting voting situations via~$N$. 
If we substitute $d$~variables $(n_1,\ldots , n_d )$ 
by $D$ new variables $(N_1,\ldots , N_D)$, say by 
setting $N_{i}$ to be the sum of $k_i$ of the variables $n_j$, for $i=1,\ldots, D$,
then we count for each $D$-tuple 
\begin{equation} \label{eq:polynomial}
p(N_1,\ldots , N_{D}) = \prod_{i=1}^{D} \binom{N_i + k_i - 1}{k_i-1}
\end{equation}
many voting situations.

In the example above, with four candidates and candidate $\cand{a}$ being the
Condorcet winner, we have $d=24$, $D=8$ and we obtain a degree~$16$
polynomial 
$$
\binom{n_{\mathsf{a}} + 5}{5}
(n_{\mathsf{ba}}+1) (n_{\mathsf{ca}}+1) (n_{\mathsf{da}}+1)  (n_{\mathsf{\ast ab}}+1) (n_{\mathsf{\ast ac}}+1) (n_{\mathsf{\ast ad}}+1)
\binom{n_{\mathsf{\ast a}} + 5}{5}
$$
to count voting situations for each $8$-tuple
$$
\left(
n_{\mathsf{a}}, n_{\mathsf{ba}}, n_{\mathsf{ca}},
n_{\mathsf{da}}, n_{\mathsf{\ast ab}}, n_{\mathsf{\ast ac}}, 
n_{\mathsf{\ast ad}}, n_{\mathsf{\ast a}} 
\right)
.
$$

Geometrically, 
the polyhedral cone $\MP\subset \R^d$ 
is replaced by a new polyhedral cone $\MP'\subset \R^{D}$ in a lower
dimension. 
As the counting is changed we obtain for the probability~\eqref{eqn:abstract-prob} of 
voting situations in $\MP$ among those in $\MS$:
\begin{equation}  \label{eqn:poly-prob-reformulation}
\Prob(n)
\; = \;
\frac{\displaystyle \sum_{x \in \MP  \cap L^d_n \cap \Z^d}
  1}{\displaystyle \sum_{x \in \MS \cap L^d_n \cap \Z^d} 1}
\; = \;
\frac{\displaystyle \sum_{y \in \MP' \cap L^D_n \cap  \Z^D}
  p(y)}{\displaystyle \sum_{y\in  \MS' \cap L^D_n \cap \Z^D} p(y)}
.
\end{equation}
Here, $\MS'$ is equal to the corresponding homogeneous polyhedral cone 
obtained from $\MS\subset\R^d$,
and $p(y)$ is the polynomial~\eqref{eq:polynomial} in $D$~variables.
In the example of Condorcet's paradox, $\MS'$ is simply equal to the
full orthant $\R^D_{\geq 0}$.

As seen in Section~\ref{sec:likelihood}, we can use Ehrhart's theory to determine an explicit formula for $\Prob(n)$. 
The right hand side of the
formula above suggests that we can do this also via {\em weighted lattice
  point counting} in dimension~$D$. A corresponding Ehrhart-type
theory has recently been considered
(see~\cite{bbdmv-2010}). 
A first implementation is available in the
package {\tt barvinok} via the command {\tt barvinok\_summate}.
We successfully tested the software on some reformulations of
three candidate elections, but so far {\tt barvinok} seems not capable to do
computations for the four candidate case. 
However, there still seems 
quite some improvement possible in the current implementation
(personal communication with Sven Verdoolaege).
It is expected that future versions of {\tt LattE integrale} will be capable of
these computations (personal communication with Matthias K\"oppe).

\subsection*{Limiting probabilities via integration}

If we want to compute the exact value of 
$\lim_{n\to \infty} \Prob(n)$ we may use integration.
Using~\eqref{eqn:poly-prob-reformulation} we get
through substitution of $y=n z$:

\begin{eqnarray*}
\lim_{n\to \infty} \Prob(n)
& = &
\lim_{n\to \infty}
\frac{\displaystyle \sum_{y \in \MP' \cap L^D_n \cap  \Z^D}
  p(y)}{\displaystyle \sum_{y \in  \MS' \cap L^D_n \cap \Z^D}
  p(y)} 
 \; = \;  
\lim_{n\to \infty}
\frac{\displaystyle \sum_{z \in \MP' \cap L^D_1 \cap  (\Z/n)^D}
  p(n z)  }{\displaystyle \sum_{z \in  \MS' \cap L^D_1 \cap (\Z/n)^D}
  p(n z)  }
\\
& = &
\lim_{n\to \infty}
\frac{\displaystyle \sum_{z \in \MP' \cap L^D_1 \cap  (\Z/n)^D}
  p(n z) / n^{\deg p} }{\displaystyle \sum_{z \in  \MS' \cap L^D_1 \cap (\Z/n)^D}
  p(n z) / n^{\deg p} }
 \; = \;  
\frac{\displaystyle\int_{\MP' \cap L^D_1}  \lt(z) \; dz}{\displaystyle\int_{\MS' \cap L^D_1} \lt(z) \; dz}
.
\end{eqnarray*}
Here, the division of numerator and denominator by
a degree of~$p$ ($\deg p$) power of~$n$
shows that the integrals on the right are taken over the
leading term $\lt(z)$ of the polynomial $p(z)$ only. 
Thus determining the exact limiting probability is achieved by
integrating a degree $d-D$ monomial over a bounded 
polyhedron ({\em polytope}) 
in the $(D-1)$-dimensional affine space~$L^D_1$. 
We refer to~\cite{ddkmpw-2011} for background on efficient integration methods
(cf.~\cite{bbdmv-2011} and \cite{schechter-1998}).

As in the case of relative volume computations in dimension~$d$,
the integral is taken with respect to the relative Lebesgue measure -- here on the affine space~$L^D_1$. 
In fact, as we are computing a quotient, any measure being a multiple of
the standard Lebesgue measure on~$L^D_1$ will give the same value.

For the example with candidate $\cand{a}$ being a Condorcet winner in
a four candidate election, the leading term to be integrated is
$$
n_{\mathsf{a}}^5 \cdot n_{\mathsf{ba}} \cdot n_{\mathsf{ca}}  \cdot
n_{\mathsf{da}} \cdot
n_{\mathsf{\ast ab}} \cdot n_{\mathsf{\ast ac}} \cdot n_{\mathsf{\ast
    ad}} \cdot
n_{\mathsf{\ast a}}^5
,
$$
which is much simpler than the full polynomial. Integrating this polynomial
over the reduced $8$-dimensional polyhedron can be done using 
{\tt LattE integrale} (called with option {\tt valuation=integrate}). 
In this way one obtains in a few seconds 
an exact value  of $1717/2048$ for the probability 
that a Condorcet winner exists (as $n$ tends to infinity).
This value corresponds to the one obtained by Gehrlein in
\cite{gehrlein-2001}
and serves as a test case for our method.
The corresponding volume computation with {\tt LattE integrale}
 (called with option {\tt valuation=volume})
in $24$~variables did not finish after several weeks of computation.
Bogdan Ichim reports (November 2011) that this volume computation is doable with his software {\tt Normaliz} (see~\cite{normaliz-link}). Nevertheless, the volume computation is much slower than the corresponding integration over the $8$-dimensional polyhedron.

\medskip

In a similar way we can deal with other voting situations as well.

\subsection*{Condorcet efficiency of plurality voting}

Assuming candidate $\cand{a}$ is a Condorcet winner, but candidate
$\cand{b}$ wins a plurality voting, we obtain a reduced system 
in the three candidate case with five variables:

\begin{eqnarray*}
\phantom{- }
n_{\mathsf{a}}  - n_{\mathsf{ba}} - n_{\mathsf{bc}} - n_{\mathsf{cb}} + n_{\mathsf{ca}}  & > & 0  \\  
\phantom{- } 
n_{\mathsf{a}}  + n_{\mathsf{ba}} - n_{\mathsf{bc}} - n_{\mathsf{cb}} - n_{\mathsf{ca}}  & > & 0  \\ 
- n_{\mathsf{a}} + n_{\mathsf{ba}} + n_{\mathsf{bc}} \;\,
\phantom{  - n_{\mathsf{cb}} - n_{\mathsf{ca}} } & > & 0  \\ 
\phantom{-n_{\mathsf{a}}  + n_{\mathsf{ba}} +}
n_{\mathsf{ba}} + n_{\mathsf{bc}} - n_{\mathsf{cb}} - n_{\mathsf{ca}}  & > & 0  
\end{eqnarray*}

Here the only reduction is the grouping
$n_{\mathsf{a}}=n_{\mathsf{ab}}+n_{\mathsf{ac}}$.
The corresponding polynomial weight is  $n_{\mathsf{a}}+1$.

The four candidate case is more involved. The linear system with $24$
variables has a comparatively small symmetry group of order~$92160$. 
We can group six variables into $n_{\mathsf{a}}$. 
Taking the reduced system~\eqref{eqn:a-condorcet-winner-among-four} 
of three inequalities with $8$~variables 
(modeling that candidate~$\cand{a}$ is a Condorcet winner) 
we have to add three inequalities for the condition that candidate $\cand{b}$ wins plurality.
These can be shortly described by $n_{\mathsf{b}}>n_{\mathsf{a}}, n_{\mathsf{c}},
n_{\mathsf{d}}$, but a grouping of variables in $n_{\mathsf{b}}$,
$n_{\mathsf{c}}$ and $n_{\mathsf{d}}$ is incompatible with the
other three conditions. Instead we use new variables 
$n_{\mathsf{b\ast a}}$, $n_{\mathsf{c\ast a}}$ and $n_{\mathsf{d\ast
    a}}$
(in~\eqref{eqn:a-condorcet-winner-among-four} combined in~$n_{\mathsf{\ast a}}$)
for preferences in which $\cand{a}$ is ranked last. Additionally we
have to keep the
variables where candidate~$\cand{a}$ is ranked third
(in~\eqref{eqn:a-condorcet-winner-among-four} combined
in~$n_{\mathsf{\ast ab}}$, $n_{\mathsf{\ast ac}}$, $n_{\mathsf{\ast ad}}$).

In the three inequalities~\eqref{eqn:a-condorcet-winner-among-four}
we can simply substitute $n_{\mathsf{\ast a}}$ by 
$n_{\mathsf{b\ast a}}+n_{\mathsf{c\ast a}}+n_{\mathsf{d\ast a}}$
and $n_{\mathsf{\ast ad}}$, $n_{\mathsf{\ast ac}}$ and $n_{\mathsf{\ast ab}}$
by $n_{\mathsf{bca}}+n_{\mathsf{cba}}$, 
$n_{\mathsf{bda}}+n_{\mathsf{dba}}$ and 
$n_{\mathsf{cda}} + n_{\mathsf{dca}}$.
The additional three linear inequalities for candidate $\cand{b}$
being a plurality winner are then:
\begin{eqnarray*}
n_{\mathsf{b\ast a}} + n_{\mathsf{ba}} + n_{\mathsf{bca}} + n_{\mathsf{bda}}  
\;\; - n_{\mathsf{a}} \qquad\qquad\qquad\qquad\quad\;\,
& > & 0 \\
n_{\mathsf{b\ast a}} + n_{\mathsf{ba}} + n_{\mathsf{bca}} + n_{\mathsf{bda}}  
\;\;  - n_{\mathsf{c\ast a}} - \, n_{\mathsf{ca}} - n_{\mathsf{cba}} - n_{\mathsf{cda}}  
& > & 0 \\
n_{\mathsf{b\ast a}} + n_{\mathsf{ba}} + n_{\mathsf{bca}} + n_{\mathsf{bda}}  
\;\;  - n_{\mathsf{d\ast a}} - n_{\mathsf{da}} - n_{\mathsf{dba}} - n_{\mathsf{dca}}  
& > & 0 
\end{eqnarray*}

This reduced linear system has $6$~inequalities for $13$~variables.
It still has a symmetry of order~$2$ coming from an interchangeable 
role of candidates~$\cand{c}$ and~$\cand{d}$.
The degree~$11$ polynomial used for integration is
$$
n_{\mathsf{a}}^5 \cdot n_{\mathsf{ba}} \cdot n_{\mathsf{ca}}  \cdot
n_{\mathsf{da}} \cdot n_{\mathsf{b\ast a}} \cdot n_{\mathsf{c\ast a}}
\cdot n_{\mathsf{d\ast a}}
.
$$
With it, using {\tt LattE integrale},
we obtain an exact limit of 
$$
\frac{10658098255011916449318509}{14352135440302080000000000}
\; = \; 
74.261410\ldots \%   
$$
for the Condorcet efficiency of plurality voting with four candidates.
To the best of our knowledge this value has not been computed before.

\subsection*{Plurality vs Plurality Runoff}

The case of Plurality vs Plurality Runoff 
has a high degree of symmetry.
For three candidates we
obtain a reduced four dimensional reformulation:

\begin{eqnarray*}
n_{\mathsf{b}} - n_{\mathsf{a}} & > & 0
\\
n_{\mathsf{a}} - n_{\mathsf{ca}} - n_{\mathsf{cb}} & > & 0
\\
n_{\mathsf{a}} + n_{\mathsf{ca}}  - n_{\mathsf{b}} - n_{\mathsf{cb}} & > & 0
\end{eqnarray*}

Counting is done via the polynomial weight
$(n_{\mathsf{a}}+1)(n_{\mathsf{b}}+1)$. Integration of $n_{\mathsf{a}} n_{\mathsf{b}}$
over the corresponding $3$-dimensional polyhedron yields the known limiting probability.

If we consider elections with $m$ candidates, $m\geq 4$, 
we can set up a
linear system with only $2(m-1)$ variables 
and $m$~inequalities.
We denote the candidates by 
$\cand{a}, \cand{b}$ and $\cand{c_i}$ for $i=1,\ldots, m-2$:

\begin{eqnarray*}
n_{\mathsf{b}} - n_{\mathsf{a}} 
& > & 0
\\
\mbox{For } i=1,\ldots, m-2: \qquad
n_{\mathsf{a}} 
- n_{\mathsf{c_i\cdot a\cdot b}} - n_{\mathsf{c_i\cdot b \cdot a}} 
& > & 0
\\
n_{\mathsf{a}} + \sum_{i=1}^{m-2} n_{\mathsf{c_i\cdot a\cdot b}} 
- n_{\mathsf{b}} - \sum_{i=1}^{m-2} n_{\mathsf{c_i\cdot b\cdot a}} 
& > & 0
\end{eqnarray*}

The first two lines model that candidate $\cand{b}$ wins plurality
over candidate~$\cand{a}$ and that candidate~$\cand{a}$ is second,
winning over candidates $\cand{c_i}$, for $i=1,\ldots,m-2$. 
The last inequality models the condition that candidate~$\cand{a}$
beats $\cand{b}$ in a pairwise comparison.
The variable $n_{\mathsf{c_i\cdot a\cdot b}}$ gives
the number of voters with candidate~$\cand{c_i}$ being their first preference and  
candidate~$\cand{a}$ being ranked before candidate~$\cand{b}$.
Similarly, $n_{\mathsf{c_i\cdot b\cdot a}}$ is the number of voters
with first preference $\cand{c_i}$  and 
candidate~$\cand{b}$ being ranked before candidate~$\cand{a}$.
We use ``$\cdot$'' to denote any ordering of candidates;
in contrast to ``$\ast$'' used before we also allow 
an empty list here.
For both variables, $n_{\mathsf{c_i\cdot a\cdot b}}$ and $n_{\mathsf{c_i\cdot b\cdot a}}$,
we group $(m-1)!/2$ of the $m!$ former variables.
The new variables $n_{\mathsf{a}}$ and $n_{\mathsf{b}}$ both represent
$(m-1)!$ former variables. Therefore, counting is
adapted using the polynomial weight
$$
\left( n_{\mathsf{a}} \cdot n_{\mathsf{b}}\right)^{(m-1)!-1} \cdot 
\prod_{i=1}^{m-2} \left( n_{\mathsf{c_i\cdot a\cdot b}} \cdot n_{\mathsf{c_i\cdot b\cdot a}}\right)^{(m-1)!/2-1} 
$$
of degree $m!-2m+2$.

The above inequalities assume that candidates~$\cand{b}$ and
$\cand{a}$ are ranked first and second in a plurality voting. So
having the probability for the corresponding voting situations,
we have to multiply by $m(m-1)$ to get the overall probability 
of a plurality winner losing in a second Plurality Runoff round.

For four candidates ($m=4$) we obtain an exact limiting probability of
$$
\frac{2988379676768359}{12173449145352192} \; = \; 24.548339\ldots \%
.
$$
This result can be obtained using the weighted, dimension-reduced
problem with {\tt LattE integrale}, as well as 
by a relative volume computation in $24$~variables. However, the
latter is a few hundred times slower than integration over
the dimension reduced polyhedron.
A similar result from a volume computation is obtained in~\cite{ddkmpw-2011}.

To be certain about our new results, we computed the
value above, as well as the likelihood for the existence of a Condorcet winner,
with a fully independent {\tt Maple} calculation, using the package {\tt Convex} (see~\cite{convex-link}). For it, we 
first obtained a {\em triangulation} (non-overlapping union of {\em simplices}) 
of the dimension-reduced polyhedron
and then applied symbolic integration to each simplex.

We also tried to solve the five candidate case, where the polyhedron
is only $7$-dimensional (in $8$~variables). 
The integration of a polynomial of degree~$112$,
however, seems a bit too difficult 
for the currently available technology.
Nevertheless it seems that we are close to obtain exact five candidate
results as well.

\section{Conclusions}

Using symmetry of linear systems we can obtain
symmetry reduced lower dimensional reformulations. 
These allow to compute exact limiting probabilities for large elections with four candidates.
In this work we only gave a few starting examples. Similar calculations are
possible for many other voting situations as well.

For the lower-dimensional weighted lattice point problems, 
efficient mathematical software for the computation of Ehrhart
quasi-polynomials will soon be available.
We anticipate that it will allow to obtain explicit formulas 
for the probability of certain voting outcomes with
four candidates and any number of voters. 
Such formulas will most likely be quite huge and 
hardly usable without computer assistance.

For elections with five or more candidates further ideas seem
necessary. One possibility to reduce the complexity of computations
further is the use of additional symmetries which remain in our
reduced systems.

\section*{Acknowledgement}

The author supervised two Bachelor projects \cite{tabak-2010} and
\cite{schreuders-2011} at Delft University of Technology on topics related to this article.
He likes to thank their two authors Frank Tabak and Marijn Schreuders 
as well as Thomas Rehn for help with some of the computations.
The author likes to thank the two anonymous referees,
as well as William~V.~Gehrlein, Karen Aardal, Jesus De Loera, Matthias K{\"o}ppe, Bogdan Ichim, 
Lars Schewe and Sven Verdoolaege for their helpful comments.

%
%

\newcommand{\etalchar}[1]{$^{#1}$}
\providecommand{\bysame}{\leavevmode\hbox to3em{\hrulefill}\thinspace}
\providecommand{\MR}{\relax\ifhmode\unskip\space\fi MR }
\providecommand{\MRhref}[2]{%
  \href{http://www.ams.org/mathscinet-getitem?mr=#1}{#2}
}
\providecommand{\href}[2]{#2}

\renewcommand\refname{Software}

\end{document}